\documentclass [10pt]{amsart} 
\usepackage {amsmath,amssymb, amscd, amsthm, epsfig, color}

\input xy
\xyoption {all}

\def \bb {\mathbb}
\def \s {{\star}}
\def \n {{\mathcal N}}

\def \endv{{\text {End }\mathcal V}}
\def \gbar {{\bar g}}
\def \ng {{\mathcal N_g}}
\def \mg {{\mathcal M_g}}
\def \p {{\mathbb P}}

\newtheorem{theorem}{Theorem} 
\newtheorem {lemma}{Lemma}

\theoremstyle{definition}

\theoremstyle {definition}

\begin{document}

\title[Intersections on the moduli space of rank two bundles]{On the intersection theory of the moduli 
space of rank two bundles}
\author {Alina Marian}
\address {Department of Mathematics}
\address {Yale University}
\email {alina.marian@yale.edu}
\author {Dragos Oprea}
\address {Department of Mathematics}
\address {Massachusetts Institute of Technology.}
\email {oprea@math.mit.edu}
\date{} 

\begin {abstract}

We give an algebro-geometric derivation of the known intersection theory on the moduli
space of stable rank $2$ bundles of odd degree over a smooth curve of genus $g$. We 
lift the computation from the moduli space to a 
Quot scheme, where we obtain the intersections by equivariant localization 
with respect to a natural torus action.

\end {abstract}

\maketitle

We compute the intersection numbers on the moduli space of stable rank $2$
odd degree bundles over a smooth curve of genus $g$. This problem has been
intensely studied in the physics and mathematics literature. Complete and
mathematically rigorous answers were obtained by Thaddeus \cite
{thaddeus}, Donaldson \cite {donaldson}, Zagier \cite {zagier}, and others
in rank $2$, and by Jeffrey-Kirwan \cite {jeffreykirwan} in arbitrary
rank. These answers are in agreement with the formulas written down by
Witten \cite {witten}. We indicate yet another method of calculation which
recovers the exact formulas obtained by the aforementioned authors in rank
$2$. Our approach works in principle in any rank, and we will turn to this
general case in future work.

To set the stage, we let $C$ be a smooth complex algebraic curve of genus $g$. We let $\mathcal N_g$ denote the moduli space of stable rank $2$ bundles with fixed odd determinant, and we write $\mathcal V$ for the universal bundle on $C \times \mathcal N_g$. We also fix, once and for all, a symplectic basis $\{1, \delta_1, \ldots, \delta_{2g}, \omega\}$ for the cohomology $H^{\star}(C)$. A result of Newstead shows that the K\"{u}nneth components of $c_2(\endv)$ generate the cohomology ring $H^{\star}(\mathcal N_g)$ \cite {newstead}. We will use the notation of Newstead and Thaddeus, writing $$c_2(\endv)=- \beta\otimes 1 +4 \sum_{k=1}^{2g} \psi_k \otimes \delta_k + 2\alpha \otimes \omega$$ for classes $\alpha\in H^2(\ng), \beta\in H^{4}(\ng), \psi_k \in H^3(\ng).$ Thaddeus showed that nonzero top intersections on $\mathcal N_g$ must contain the $\psi_k$s in pairs, which can then be removed using the formula $$\int_{\ng} \alpha^m \beta^n \prod_{k=1}^{g} (\psi_k \psi_{k+g})^{p_k}=\int_{\mathcal N_{g-p}} \alpha^m \beta^n$$ where $p=\sum_k p_k$ and $2m+4n+6p=6g-6$. The top intersections of $\alpha$ and $\beta$ are further determined: \begin {theorem}\cite {thaddeus}\cite{donaldson}\cite{zagier}\cite{jeffreykirwan} \begin {equation} \label{main} \int_{\ng}\alpha^m \beta^n =(-1)^g 2^{2g-2} \frac{m!}{(m-g+1)!} (2^{m-g+1}-2)B_{m-g+1}.\end {equation} \end {theorem} Here $B_k$ are the Bernoulli numbers defined, for instance, by the power series expansion \begin{equation}\label{bernoulli} -\frac{u}{\sinh u}=\sum_{k} \frac{2^{k}-2}{k!} B_{k} u^{k}.\end{equation}

In this paper, we reprove theorem $\ref{main}$. The idea is to lift the computation from $\ng$ to a Quot scheme as indicated in the diagram below. Then, we effectively calculate the needed intersections on Quot by equivariant localization with respect to a natural torus action. The convenience of this approach lies in that the fixed loci are easy to understand; in any rank they are essentially symmetric products of $C$. In rank $2$, their total contribution can be evaluated to the intersection numbers in $\eqref{main}$. 

\begin{center}
$
\xymatrix {Quot_{N, d} \ar@{<.>}[r] & {\p_{N, d}} \ar[d]^{\pi} \\
 & \mg & \ng \times J \ar[l]^{\tau}}
$
\end{center}

The idea that the intersection theory of the moduli space of rank $r$ bundles on a
curve and that of a suitable Quot scheme are related goes back to Witten \cite 
{witten2} in the
context of the Verlinde formula. Moreover, in \cite {bertram} the authors have the
reverse approach of calculating certain intersection numbers on Quot in low rank and
genus by using the intersection theory of the moduli space of bundles. In this note
however, the translation of the intersection theory of the moduli space of bundles
into that of the Quot scheme is really very straightforward.

The spaces in the diagram are as follows.
\begin{itemize}
\item $\mg$ denotes the moduli space of rank 2, odd degree stable bundles on $C$. By 
contrast with $\ng,$ the determinant of the bundles is allowed to vary in $\mg.$ There is a finite covering 
map $$\tau: \n_g \times J \rightarrow \mg$$ of degree $4^g$, given by tensoring bundles. Here $J$ is the Jacobian of degree 0 line bundles on $C$. We will write $\widetilde {\mathcal V}$ for the universal sheaf on $\mg \times C$, which is only defined up to twisting with line bundles from $\mg$. \\

\item For a positive integer $N$, we let $\p_{N, d}$ be the projective bundle $$\pi: 
{\bb P}_{N,d} \rightarrow \mg$$ whose fiber over a stable $[V]\in \mg$ is ${{\bb P} 
(H^0 (C, V)^{\oplus N}}).$ In other words,$$\p_{N, d}  = {\bb P} (\text {pr}_{\s} 
\widetilde {\mathcal V}^{\oplus N} ).$$ 
We will take the degree to be as large as convenient to ensure, for instance, that the fiber dimension of $\p_{N, d}$ is constant. One can view $\p_{N, d}$ as a fine moduli space of pairs $(V, {\phi})$ consisting of a stable rank 2, degree $d$ bundle $V$ together with a nonzero $N$-tuple of holomorphic sections $${\phi} = (\phi^1, \ldots \phi^N):{\mathcal {O}}^N \to V$$ considered projectively. As such, there is a universal morphism $$\Phi: {\mathcal {O}}^N\to \overline {\mathcal V} \text{ on }\p_{N, d}\times C.$$ 

A whole series of moduli spaces of pairs of vector bundles with $N$ sections, each labeled by a parameter $\tau$, were defined and studied in \cite{thaddeus2} for $N=1$ and \cite{bertram} for arbitrary 
$N$. $\p_{N, d}$ is one of these moduli spaces, undoubtedly the least exciting one from the point of view of 
studying a new object, due to its straightforward relationship with $\mg$. 

One would expect that the universal bundle $\overline {\mathcal V}$ on the total space of $\p_{N,d}$ and the 
noncanonical universal bundle $\widetilde {\mathcal V}$ on its base $\mg$ to be 
closely related, and indeed they 
coincide up to a twist of ${\mathcal {O}}(1)$ \begin{equation}\label{universal}\overline {\mathcal 
V}=\pi^{\star} \widetilde {\mathcal V}\otimes {\mathcal {O}}_{\p_{N,d}}(1).\end{equation}

\item Finally, $Quot_{N, d}$ is Grothendieck's Quot scheme parametrizing degree $d$
rank $2$ subbundles $E$ of the trivial rank $N$ bundle on $C$, $$0 \to E \to {\mathcal {O}}^N\to
F\to 0.$$ We let $$0\to \mathcal E\to {\mathcal {O}}^N\to \mathcal F\to 0$$ be the universal
sequence on $Quot_{N, d}\times C$. 

For large $d$ relative to $N$ and $g$, the Quot
scheme is irreducible, generically smooth of the expected dimension 
$Nd-2(N-2)(g-1)$ \cite {bertram}. Then, $Quot_{N, d}$ and
$\p_{N, d}$ are birational and they agree on the open subscheme corresponding to
subbundles $E=V^{\vee}$ where $\phi: {\mathcal {O}}^N \to V$ is generically surjective. The
universal structures also coincide on this open set. For arbitrary $d$, $Quot_{N,d}$
may be badly behaved, but intersection numbers can be defined with the aid of the
virtual fundamental cycle constructed in \cite {opreamarian}.

\end {itemize}

We now define the cohomology classes that we are going to intersect. We consider the
K\"{u}nneth decomposition $$c_2 (\text {End} \mathcal {\widetilde V}) = -
\tilde{\beta} \otimes 1 + 4 \sum_{k=1}^{2g} \tilde \psi_k \otimes \delta_{k} + 2
\tilde{\alpha} \otimes \omega$$ of the universal endomorphism bundle on $\mg \times
C.$ In keeping with the notation of \cite{jeffreykirwan}, we further let
$$c_i(\widetilde {\mathcal V}) = \tilde a_i \otimes 1 + \sum_{k=1}^{2g} {\tilde b_i^k}
\otimes \delta_k + \tilde f_i \otimes \omega, \, 1\leq i \leq 2 $$ be the K\"unneth 
decomposition of the
(noncanonical) universal bundle $\widetilde {\mathcal V}$ on $\mg \times C.$ Then
$$\tilde f_1 = d, \, \, \tilde{\beta} = \tilde a_1^2 - 4 \tilde a_2, \, \,
\tilde{\alpha} = 2 \tilde f_2 + \sum_{k=1}^{g} \tilde b_1^k \tilde b_1^{k+g} - d
\tilde a_1.$$ It is an easy exercise, using that $\ng$ is simply connected, to see
that \begin{equation}\label{pullback}\tau^{\star}\tilde {\alpha}=\alpha,\;\;\;
\tau^{\star}\tilde {\beta}=\beta,\;\;\; \tau^{\star} (\sum_{k=1}^{g} \tilde b_1^k
\tilde b_1^{k+g}) = 4 \theta,\end{equation} where $\theta$ is the class of the theta
divisor on $J.$

On $\p_{N, d}$, we let $\zeta$ denote the first Chern class of ${\mathcal {O}}(1)$ and let
\begin{equation}c_i(\overline {\mathcal V}) = {\bar{a}}_i \otimes 1 + \sum_{k=1}^{2g} {\bar{b}}_i^k 
\otimes \delta_k + {\bar{f}}_i
\otimes
\omega, \, 1 \leq i \leq 2 \end{equation} be the K{\"{u}}nneth decomposition of the 
Chern classes of $\overline {\mathcal V}.$ 

Finally, we consider the corresponding $a, b, f$ classes on $Quot_{N, d}$: $$c_i(\mathcal E^{\vee})=a_i\otimes 1+ \sum_{k=1}^{2g} b_i^{k}\otimes \delta_k + f_i \otimes \omega, \, 1\leq i\leq 2.$$

We now show how to express any top intersection of $\alpha$ and $\beta$ classes on $\ng$ as an intersection on $Quot_{N, d}$. Let $m$ and $n$ be any nonnegative integers  
such that $$m+2n = 4g-3, \,\,\, m\geq g.$$ With the aid of \eqref{pullback} we note 
\begin{equation}\label{red1}\int_{\mg} (\tilde{\alpha} + \sum_{k=1}^{g} \tilde b_1^k \tilde b_1^{k+g})^m \tilde{\beta}^n = 
\frac{1}{4^g} \cdot \binom{m}{g} \int_{J} (4 \theta)^g  \int_{\n_g} 
\alpha^{m-g} \beta^n.\end {equation}
 
A top intersection on $\mg$ which is invariant under the normalization of the 
universal bundle $\mathcal {\widetilde V}$ on $\mg \times C$ can be readily expressed 
as a
top intersection on $\p_{N, d}$ as follows. We assume from now on that $N$ is odd and we let $$2M = N (d-2\bar{g})-1.$$
Then
\begin {equation}\label{red2}\int_{\mg} \left(\tilde{\alpha} + \sum_{k=1}^{g} \tilde b_1^k \tilde b_1^{k+g}\right)^m \tilde{\beta}^n = 
\int_{\mg} \left( 2\tilde f_2 + 2 \sum_{k=1}^{g} \tilde b_1^k \tilde b_1^{k+g}-d \tilde a_1 \right)^m (\tilde a_1^2-4\tilde a_2)^n=\end {equation}

$$=\int_{\p_{N, d}} \zeta^{2M} \left( 2{\bar{f}}_2 + 2 \sum_{k=1}^{g} 
{\bar{b}}_1^k 
{\bar{b}}_1^{k+g}-d {\bar{a}}_1 \right)^m \left({\bar{a}}_1^2-4{\bar{a}}_2\right)^n =$$ $$=\int_{\p_{N, d}} {\bar{a}}_2^M 
\left ( 2{\bar{f}}_2 + 2 \sum_{k=1}^{g}
{\bar{b}}_1^k
{\bar{b}}_1^{k+g} - d {\bar{a}}_1 \right)^m \left({\bar{a}}_1^2-4{\bar{a}}_2\right)^n.$$
Here we used \eqref{universal} to write $${\bar{a}}_2^M = (\zeta^2 + \tilde a_1 \zeta + \tilde a_2)^M
= \zeta^{2M} + \text { lower order terms in }\zeta,$$ the latter being zero when paired with a top
intersection from the base $\mg$. 

Finally, this last intersection number can be
transfered to $Quot_{N, d}$ using the results of \cite {marian}. It is shown there that the equality \begin{equation}\label{red3} \int_{\p_{N, d}} {\bar{a}}_2^M R(\bar a, \bar b, \bar
f)=\int_{Quot_{N, d}} a_2^M R(a, b, f)\end {equation} holds for any polynomial $R$ in the $a, b, f$ classes, in the regime that $N$ is large compared to the genus $g$, and in turn $d$ is large enough relative to $N$ and $g$, so that $Quot_{N,d}$ is irreducible of the expected dimension. Moreover, the equality \begin{equation}\label{virdeg}\int_{\left[Quot_{N,d}\right]^{vir}}a_2^M R(a, b, f)=\int_{\left[Quot_{N, d-2}\right]^{vir}} a_2^{M-N} R(a, b, f)\end {equation} established in \cite {opreamarian} allows us to assume from this moment on that the degree $d$ is as small as desired relative to $N$. The trade-off is that we need to make use of the virtual fundamental classes alluded to above and defined in \cite {opreamarian}. 

Putting $\eqref{red1}$,
$\eqref{red2}$, $\eqref{red3}$ together we obtain \begin {equation} \int_{\ng}
\alpha^{m-g}\beta^n=\frac{(m-g)!}{m!}\int_{\left[Quot_{N, d}\right]^{vir}} a_2^M \left (2f_2 + 2 \sum_{k=1}^{g}b_1^k b_1^{k+g}-d {a}_1\right)^m \left(a_1^2-4a_2\right)^n.\end
{equation} To prove $\eqref {main}$, we will show that \begin
{equation}\label{red}\int_{[Quot_{N, d}]^{vir}} a_2^M \left (2f_2 + 2
\sum_{k=1}^{g}b_1^k b_1^{k+g}-d {a}_1\right)^m \left(a_1^2-4a_2\right)^n=\end {equation} $$=(-1)^g(2^{m-1}-2^{2g-1})
\frac{m!}{(m-2g+1)!} B_{m-2g+1}.$$

Equation $\eqref{red}$ will be verified by virtual localization. The torus action we will use
and its fixed loci were described in \cite {opreamarian}. For the convenience of the
reader, we summarize the facts we need below.

The torus action on $Quot_{N, d}$ is induced by the fiberwise $\mathbb C^{\star}$ action on ${\mathcal {O}}^N$ with distinct weights $-\lambda_1, \ldots, -\lambda_N$. On closed points, the action of $g\in \mathbb C^{\star}$ is $$\left[E\stackrel{i}{\to} {\mathcal {O}}^N \right]\mapsto
\left[E\stackrel{g\circ i}{\to}{\mathcal {O}}^N\right].$$ The fixed loci $Z$ correspond to split subbundles $$E= L_1 \oplus L_2$$ where $L_1$ and $L_2$ are line subbundles of copies of ${\mathcal {O}}$ of degree $-d_1$ and $-d_2$. Thus $Z=\text {Sym}^{d_1}
C \times \text {Sym}^{d_2} C$ and $$\mathcal E|_{Z}=\mathcal L_1\oplus \mathcal L_2$$
where we let $\mathcal L_1$ and $\mathcal L_2$ be the universal line subbundles on
$\text {Sym}^{d_1} C\times C$ and $\text {Sym}^{d_2} C\times C$. We write $$c_1(\mathcal
L_i^{\vee})=x_i\otimes 1+ \sum_{k} y_i^{k}\otimes \delta_k+d_i \otimes \omega,\,\,\, 1\leq i\leq 2.$$ 

We set the weights to be the $N^{\text {th}}$ roots of unity. The equivariant Euler class of the
virtual normal bundle of $Z$ in $Quot_{N, d}$ was determined in \cite {opreamarian} to be \begin
{equation}\label{normal} \frac{1}{e_{T}(\n^{vir})} = (-1)^{g} \left ((\lambda_1 h +x_1)- (\lambda_2
h + x_2\right) )^{-2\gbar}\cdot {\prod_{i=1}^{2}} \left ( \frac{x_i}{(\lambda_i h + x_i)^N - h^N}
\right )^{d_i - \gbar}\end {equation}$$ \cdot\prod_{i=1}^{2}\exp\left(\theta_i\cdot{\left ( \frac{N
(\lambda_i h + x_i)^{N-1}}{(\lambda_i h + x_i)^N - h^N} - \frac{1}{x_i}\right )}\right).$$ Here
$\gbar=g-1$, $h$ is the equivariant parameter, and $\theta_i$ are the 
pullbacks of the theta divisor
from the Jacobian. Moreover, it is clear that \begin{equation}\label{rest1}a_1 = (x_1 + \lambda_1 h)
+ (x_2+\lambda_2 h), \;\; a_2=(x_1+\lambda_1 h) (x_2+\lambda_2 h),\end{equation}
\begin{equation}\label{rest2}b_{1}^{k}=y_1^{k}+y_2^{k},\end{equation}
\begin{equation}\label{rest3}f_2=-\sum_{k=1}^{2g} y_1^{k} y_{2}^{k+g} + d_2(x_1+\lambda_1 h) +
d_1(x_2+\lambda_2 h).\end{equation} We collect $\eqref{normal}$, $\eqref{rest1}$, $\eqref{rest2}$,
$\eqref{rest3}$, and rewrite the left hand side of $\eqref{red}$, via the virtual localization
theorem, \begin{equation}\label{bigsum}\text {LHS of }\eqref{red}=(-1)^{g} \sum_{d_1, d_2,
\lambda_1, \lambda_2}\mathcal I_{d_1, d_2, \lambda_1, \lambda_2}\end{equation} where the sum ranges
over all degree splittings $d=d_1+d_2$ and pairs of distinct roots of unity $(\lambda_1,
\lambda_2)$. The summand $\mathcal I_{d_1, d_2, \lambda_1, \lambda_2}$ is defined as the evaluation
on $\text {Sym}^{d_1} C\times \text {Sym}^{d_2} C$ of the expression
  
\begin{equation}\label{tosum}
\left((\lambda_1 h + x_1)-(\lambda_2 h+ x_2)\right)^{2n-2\gbar}
\left(2\theta_1+2\theta_2+(d_2-d_1) ((\lambda_1 h+x_1)-(\lambda_2 h
+x_2))\right)^{m}\end {equation} $$\cdot \prod_{i=1}^{2}(\lambda_i h+ x_i)^{M} \left (
\frac{x_i}{(\lambda_i h + x_i)^N - h^N} \right )^{d_i - \gbar}
\exp\left(\theta_i\cdot{\left ( \frac{N (\lambda_i h + x_i)^{N-1}}{(\lambda_i h + x_i)^N
- h^N} - \frac{1}{x_i}\right )}\right).$$
\vspace{.12in}

To carry out this evaluation, we use the following standard facts regarding intersections on $\text {Sym}^{d}C$:
\begin{equation} \label{xtheta} x^{d-l} \theta^{l} = \frac{g!}{(g-l)!}\, \,
\text{for} \, \, l \leq g, \, \, \text{and} \, \, x^{d-l} \theta^{l} = 0 \text{ for}
\, \, l > g.  \end{equation}
We will henceforth replace any $\theta^{l}$ appearing in a top intersection on $\text {Sym}^d C$ by
$\frac{g!}{(g-l)!}x^{l}$. Then

\begin {equation}\label{thetay} \frac
{\theta^{l}}{l!} \exp\left(\theta\cdot{\left ( \frac{N (\lambda h +
x)^{N-1}}{(\lambda h + x)^N - h^N} - \frac{1}{x}\right )}\right)= \sum_{k}
\frac{\theta^{l+k}}{l!\, k!} \left( \frac{N (\lambda h + x)^{N-1}}{(\lambda h +
x)^N - h^N} - \frac{1}{x}\right )^{k}= \end{equation} $$=\sum_{k\leq g-l} \frac{g!\,
x^{l+k}}{l!\,(g-l-k)!\,k!} \left( \frac{N (\lambda h + x)^{N-1}}{(\lambda h + x)^N
- h^N} - \frac{1}{x}\right )^{k}= N^{g-l} \binom{g}{l} x^{g}\cdot
\frac{(\lambda h+x)^{(N-1)(g-l)}}{((\lambda h+x)^{N}-h^N)^{g-l}}.$$
\vspace{.12in}

We further set $$\bar{x}_i = \frac{x_i}{\lambda_i h} \text{ and } \mathcal J=\lambda_1 (1+\bar {x}_1)-\lambda_2(1+\bar{x}_2).$$ In terms of the rescaled variables,
${\mathcal I}_{d_1, d_2, \lambda_1, \lambda_2}$ becomes, via \eqref{thetay}, the residue at $\bar{x}_1=\bar{x}_2=0$ of
 
$$\sum_{l_1+ l_2 + s = m} \left ( 
\frac{m!}{s!} (d_2 - d_1)^s \cdot \mathcal J^{2n-2\gbar +s} \cdot \prod_{i=1}^{2} 2^{l_i}N^{g-l_i}
\binom{g}{l_i} \lambda_i^{M+1 + l_i} \frac{(1+\bar{x}_i)^{(N-1)(g-l_i) +M}}{((1+\bar{x}_i)^N 
-1)^{d_i - 
l_i +1}} \right ).$$\vspace{.1in}
We expand $$\mathcal J^{2n-2\gbar+s}=\sum_{k=0}^{d}\sum_{\alpha_1+\alpha_2=k} (-1)^{\alpha_2} \binom{k}{\alpha_1} \lambda_1^{\alpha_1} \lambda_2^{\alpha_2} \cdot {\mathfrak z}_N(s, k) \cdot \bar {x}_1^{\alpha_1} \bar{x}_2^{\alpha_2}$$ where we set $${\mathfrak z}_N(s, k)=\binom{2n-2\gbar+s}{k} (\lambda_1-\lambda_2)^{2n-2\gbar+s-k}.$$ 

In \cite {opreamarian}, we computed the residue $$\text {Res}_{x=0}\left\{ \frac{x^a (1+x)^{N-1+b}}{((1+x)^N-1)^{c+1}}\right\}=\frac{1}{N}\sum_{p=0}^{a} (-1)^{a-p}\binom{\frac{b+p}{N}}{c}\binom{a}{p}.$$ Therefore, $\mathcal I_{d_1, d_2, \lambda_1, \lambda_2}$ becomes \begin{equation}\label{tosum2}\sum_{l_1+l_2+s=m}\sum_{k=0}^{d}\sum_{\alpha_1+ \alpha_2=k}\sum_{p_1, p_2} \left\{\frac{m!}{s!}(d_2-d_1)^s (-1)^{\alpha_2}\binom{k}{\alpha_1}\cdot {\mathfrak z}_N(s, k)\cdot \right.\end{equation} $$\left.\cdot \prod_{i=1}^{2} 2^{l_i} N^{\gbar-l_i}\binom{g}{l_i}\lambda_i^{M+l_i+\alpha_i+1}\binom{\frac{M+(N-1)(\gbar-l_i)+p_i}{N}}{d_i-l_i}\binom{\alpha_i}{p_i}(-1)^{\alpha_i-p_i}\right \}.$$ 

Now $\eqref{bigsum}$ computes an intersection number on $\mg$, hence it is
independent of $N$. Thus it is enough to make $N\to \infty$ in the above
expression. Lemma $\ref{bern}$ below clarifies the $N$ dependence of the
sum over the roots of unity \begin{equation}\label{sumoverroots}\sum
{\mathfrak z_{N}}(s, k) \cdot
\lambda_1^{M+l_1+\alpha_1+1}\lambda_2^{M+l_2+\alpha_2+1},\end{equation} by
writing $\lambda_1=\zeta\lambda_2$. Lemma $\ref{binoms}$ takes care of all
the other terms in \eqref{tosum2}.

A moment's thought shows that in the limit the sum over the roots of unity
of the terms \eqref{tosum2} reduces to $$\frac{1}{2}\sum
\frac{m!}{s!}(d_2-d_1)^{s} \cdot \prod_{i=1}^{2} 2^{l_i}\binom{g}{l_i}
\cdot k! \text { Coeff }_{x_1^{\alpha_1}}
\binom{x_1+\frac{d}{2}-l_1}{d_1-l_1}\cdot \text { Coeff
}_{(-x_2)^{\alpha_2}}\binom{x_2+\frac{d}{2}-l_2}{d_2-l_2}\cdot {\mathfrak
z}_{\infty}(s, k)$$ with $${\mathfrak z}_{\infty}(s,
k)= \binom{2n-2\gbar+s}{k}\cdot
\frac{B_{2\gbar-2n-s+k}}{(2\gbar-2n-s+k)!}\cdot (1-2^{2n-2\gbar+s-k+1}).$$\vspace{.12in}
Summing over $\alpha_1+\alpha_2=k$ first, we obtain $$\frac{1}{2}\sum
\frac{m!}{s!}(d_2-d_1)^{s} \cdot \prod_{i=1}^{2} 2^{l_i}\binom{g}{l_i}
\cdot k! \text { Coeff }_{x^{k}} \binom{x +
\frac{d}{2}-l_1}{d_1-l_1}\binom{-x+\frac{d}{2}-l_2}{d_2-l_2} \cdot
{\mathfrak z}_{\infty}(s, k).$$ \vspace{.1in} We sum next over
$d_1+d_2=d$. An easy induction on $m$ shows that generally
$$\sum_{b_1+b_2=a} \binom{a_1}{b_1} \binom{a_2}{b_2}
(t+b_1-b_2)^m=(t+a_1-a_2)^m$$ holds whenever $a_1+a_2=a$. Via this
observation, our expression simplifies to $$\frac{1}{2} \sum 
\frac{m!}{s!}\cdot
\prod_{i=1}^{2} 2^{l_i}\binom{g}{l_i}\cdot k! \text { Coeff}_{x^k} (-2x)^s
\cdot {\mathfrak z}_{\infty}(s, k)= 2^{m-1} m! \sum
\binom{g}{l_1}\binom{g}{l_2} (-1)^{s}\cdot {\mathfrak z}_{\infty}(s, s)$$
$$=2^{m-1} m!\, 
\frac{B_{2\gbar-2n}}{(2\gbar-2n)!}\,(1-2^{2n-2\gbar+1})\sum
(-1)^{s} \binom{g}{l_1}\binom{g}{l_2}\binom{2n-2\gbar+s}{s}$$
$$=2^{m-1}\,\frac{m!}{(m-2g+1)!}\,B_{m-2g+1}\,(1-2^{-m+2g}) \sum
\binom{g}{l_1} \binom{g}{l_2} \binom{m-2g}{s}$$
$$=2^{m-1}\frac{m!}{(m-2g+1)!}\,B_{m-2g+1}\,(1-2^{-m+2g}).$$ \vspace{.1in}
This last equality and \eqref{bigsum} complete the proof of \eqref{red},
hence of the theorem.

\begin {lemma}\label{bern} For all integers $a$ and $k$ we have \begin {equation}\label{limit}
\lim_{N\to \infty} \frac{1}{N^{k}} \left(\sum_{\zeta\neq 1}
\frac{\zeta^{\frac{N-1}{2}+a}}{(1-\zeta)^{k}}\right)=(1-2^{-k+1})\cdot\frac{B_{k}}{k!}\end
{equation} the sum being taken over the $N^{\text {th}}$ roots $\zeta$ of $1$. \end {lemma}

{\bf Proof.} When $N$ is large, the sum to compute is $0$ for $k<0$, and $-1$ for $k=0$. We may thus assume that $k\geq 1$. We introduce the auxiliary variable $z$, and evaluate $$\sum_{k=1}^{\infty} {z^{k-1}}
\left(\sum_{\zeta\neq 1} \frac{\zeta^{\frac{N-1}{2}+a}}{(1-\zeta)^{k}}\right)=\sum_{\zeta\neq 1}
\frac{\zeta^{\frac{N-1}{2}+a}}{1-z-\zeta}=\frac{1}{z}+N
\frac{(1-z)^{\frac{N-1}{2}+a-1}}{(1-z)^{N}-1}.$$ Setting $z=\frac{u}{N}$ and making $N\to \infty$
we obtain $$\sum_{k=1}^{\infty} \frac{u^{k}}{N^k} \left(\sum_{\zeta\neq 1}
\frac{\zeta^{\frac{N-1}{2}+a}}{(1-\zeta)^{k}}\right)=1+u\cdot
\frac{(1-u/N)^{\frac{N-1}{2}+a-1}}{(1-u/N)^{N}-1}\to 1+
\frac{ue^{-u/2}}{e^{-u}-1}=1-\frac{u}{2\sinh\frac{u}{2}}.$$ The lemma follows.

\begin{lemma}\label{binoms}
Let $b, \alpha$ be fixed non-negative integers and $z$ a real number. Then the following limit $$\lim_{N\to \infty} N^{\alpha} \sum_{p=0}^{\alpha} \binom{z+\frac{p}{N}}{b}\binom{\alpha}{p}(-1)^{\alpha-p}$$ equals the coefficient of $x^{\alpha}$ in $\alpha!\, \binom{x+z}{b}.$
\end {lemma}

{\bf Proof.} Let us write $$f_b(z)=\sum_{p=0}^{\alpha} \binom{z+\frac{p}{N}}{b}\binom{\alpha}{p}(-1)^{\alpha-p}.$$ The recursion $$f_{b+1}(z+1)=f_{b+1}(z)+f_{b}(z)$$ implies by induction that $$f_b(z)=\sum_{j=0}^{b} \binom{z}{b-j}f_j(0).$$ We evaluate $$N^{\alpha}f_{j}(0)=N^{\alpha}\sum_{p=0}^{\alpha} \binom{\frac{p}{N}}{j}\binom{\alpha}{p}(-1)^{\alpha-p}=N^{\alpha}\sum_{i=0}^{j} c(j,i)\sum_{p=0}^{\alpha}\left(\frac{p}{N}\right)^{i}\binom{\alpha}{p}(-1)^{\alpha-p},$$ where $c(j,i)$ are the coefficients defined by the expansion $$\binom{x}{j}=\sum_{i=0}^{j} c(j, i)\, x^{i}.$$ We use the Euler identities \begin{equation*}\sum_{p=0}^{\alpha} p^i \binom{\alpha}{p}(-1)^{\alpha-p}=\begin{cases}0 & \text {if } i<\alpha\\ \alpha! &\text {if } i=\alpha  \end{cases}.\end{equation*} Making $N\to \infty$ we obtain $$\lim_{N\to \infty} N^{\alpha} f_j(0)=\alpha!\, c(j, \alpha).$$ Therefore, $$\lim_{N\to \infty} f_b(z)=\alpha! \,\sum_{j=0}^{b}\binom{z}{b-j} c(j, \alpha).$$ This expression can be computed as the coefficient of $ x^{\alpha}$ in $$\alpha!\,\sum_{j=0}^{b} \binom{z}{b-j}\binom{x}{j}=\alpha!\binom{x+z}{b}.$$

\end{document}